\documentclass[11pt]{amsart}
\usepackage{amsmath,amsthm,amsfonts,amssymb,amscd,hyperref}
\usepackage{fancyhdr}
\usepackage{setspace}
\newtheorem{problem}{PROBLEM}
\newtheorem{theorem}{THEOREM}[section]
\newtheorem{corollary}{Corollary}[theorem]
\newtheorem{definition}[theorem]{Definition}
\newtheorem{remark}{Remark}
\newtheorem{example}{Example}[section]
\newtheorem{lemma}[theorem]{Lemma}
\newtheorem{proposition}[theorem]{Proposition}
\numberwithin{equation}{section}
\newcommand{\C}{\mathbb{C}}

\newcommand{\x}{\zeta}
\newcommand{\y}{\eta}

\newcommand{\z}{\mu}

\newcommand{\fs}{\mathfrak{S}}

\renewcommand{\a}{\alpha}
\renewcommand{\b}{\beta}
\renewcommand{\r}{\gamma}
\newcommand{\s}{\sigma}

\newcommand{\U}{\textsl{U}}

\newcommand{\bt}{\begin{theorem}}
\newcommand{\et}{\end{theorem}}
\newcommand{\bco}{\begin{corollary}}
\newcommand{\eco}{\end{corollary}}
\newcommand{\bd}{\begin{definition}}
\newcommand{\ed}{\end{definition}}
\newcommand{\bp}{\begin{problem}}
\newcommand{\ep}{\end{problem}}
\newcommand{\bl}{\begin{lemma}}
\newcommand{\el}{\end{lemma}}
\newcommand{\bprop}{\begin{proposition}}
\newcommand{\eprop}{\end{proposition}}
\newcommand{\br}{\begin{remark}}
\newcommand{\er}{\end{remark}}
\newcommand{\bpf}{\begin{proof}}
\newcommand{\epf}{\end{proof}}
\title{Picard-Vessiot Extensions For Unipotent Algebraic Groups}
\author{V. Ravi Srinivasan}

\address{\noindent Department of Mathematics and Computer Science, Rutgers University, Newark, NJ 07102.}
\email{ravisri@rutgers.edu}

\begin{document}

\maketitle
\begin{abstract} Let $F$ be a  differential field of characteristic zero. In this article, we construct Picard-Vessiot extensions of $F$ whose differential Galois group is isomorphic to the full unipotent subgroup of the upper triangular group defined over the field of constants of $F$. We will also give a procedure to compute linear differential operators for our Picard-Vessiot extensions.  We do not require the condition that the field of constants be algebraically closed.  \end{abstract}

\section{Introduction}
Throughout this article, we fix a ground differential field $F$ of characteristic zero. All the differential fields considered henceforth are either differential subfields of $F$ or a differential field extension of $F$. We deal with differential fields equipped with only one derivation and we reserve the notation $'$ to denote that derivation map. The author assumes that the reader is familiar with the notion of a differential field and Picard-Vessiot theory. For precise definitions see \cite{Magid 1} and \cite{Marius-Singer}.

Let $\U(n,C)$ be the subgroup of $GL(n,C)$ of all upper triangular matrices with 1's on the diagonal. In this article, we describe a procedure to compute linear homogeneous differential equations over $F$ for the group $\U(n,C)$.  Here is the statement of our main result:

\begin{theorem}\label{Main result}

Let $F$ be a differential field and $C$ be its field of constants. Suppose that $F$ contains distinct elements $f_1,f_2,\cdots,f_n$ satisfying the following condition:

\begin{itemize}\item[{\bf(C)}]if there are elements $c_1,c_2,\cdots,c_n\in C$ and an element $f\in F$ such that $\sum^n_{i=1}c_if_i=f'$ then $c_i=0$ for all $i$.\end{itemize}

Let $\fs:=\{\x_{i,j}|1\leq i\leq n, 1\leq j\leq n+1-i\}$ be a set of $n(n+1)/2$-variables and let $E:=F(\fs)$ be the field of rational functions in these variables over the field $F$.  Let \begin{equation} g:=\begin{pmatrix}
1 & \x_{1,1} & \x_{2,1}&\cdots& \x_{n,1}\\
0 & 1 &  \x_{1,2}& \cdots& \x_{n-1,2}\\
\vdots & \vdots &  \vdots&\vdots& \vdots\\
0 & 0 &  \cdots&1 & \x_{1,n}\\
0 & 0 &  \cdots & 0& 1
\end{pmatrix},\ \  A:=\begin{pmatrix}
0 & f_1 & 0&\cdots& 0\\
0 & 0 &  f_2& \cdots& 0\\
\vdots & \vdots &  \vdots&\vdots& \vdots\\
0 & 0 &  \cdots&0 & f_n\\
0 & 0 &  \cdots & 0& 0
\end{pmatrix}.  \end{equation}
Extend the derivation of $F$ to $E$ by setting \begin{equation}\label{defn derv}g'=Ag.\end{equation}
 Then $E$ is a Picard-Vessiot extension of $F$ for the differential operator \begin{equation}\label{gen diff eqn}L(Y):=\frac{w(Y,1,\x_{1,1},\x_{2,1},\cdots,\x_{n,1})}{w(1,\x_{1,1},\x_{2,1},\cdots,\x_{n,1})}, \end{equation} $V:=$ span$_C\{\x_{1,1},\x_{2,1},\cdots,\x_{n,1}\}$ is the full set of solutions of the equation $L(Y)=0$, and the differential Galois group $G:\label{diff eq}=G(E|F)$ is naturally isomorphic to the group $\U(n+1,C)$.

\end{theorem}

When $F=\C(z)$ with the usual derivation $d/dz$, a maple program to compute equation \ref{gen diff eqn}
can be found in my website \cite{Ravi.S w}.

 Under the conditions that the field of constants $C$ of $F$ be algebraically closed  and that $F$ has a non zero but finite transcendence degree over $C$, Biyalinicki-Birula \cite{Bi-Bi} has shown that every connected nilpotent algebraic group defined over $C$ can be realized as a differential Galois group. And, in \cite{mit-sing} the authors construct Picard-Vessiot extensions of $F$ for connected linear algebraic groups defined over $C$. Thus, in particular, the existence of a Picard-Vessiot extension for the group $U(n+1, C)$ is known. Our approach differs and gives a different perspective from the above mentioned articles.
 


{\bf Acknowledgements.} The author would like to thank  the participants of the Kolchin Seminar, especially Phyllis Cassidy and William Sit, for their suggestions and criticisms on this work.  Special thanks to Jonathan Sparling  and to William Keigher for their encouragement and support throughout this project.

\subsection{Preliminaries}

Here we will state some results that will be used often in this article.
The following theorems aid us in constructing no new constants extension by adjoining antiderivatives, see \cite{Magid 1} or \cite{Rosenlicht-Singer}.

\bt \label{new constant} Let $F$ be a differential field and let $E=F(\x)$ be a differential field extension of $F$  such that $\x'\in F$. Suppose that $E$ has a new constant. Then there is a $y\in F$ such that $y'=\x'$.\et

The next theorem characterizes the algebraic dependence of antiderivatives and it is a special case of the Kolchin-Ostrowski theorem (see \cite{j.ax}, \cite{M.Ros}).

\bt \label{transcen}Let $  E\supset  F$ be a no new constants
extension and for $i=1,2, \cdots,n,$ let $\x_i\in  E$ be
antiderivatives of $  F$. Then either $\x_i$'s are algebraically
independent over $  F$ or there is a tuple
$(c_1, \cdots,c_n)\in C^n-\{0\}$ such that
$\sum^n_{i=1}c_i\x_i$ $\in  F$.\et

A direct consequence of the above theorem is the following proposition, see \cite{Ravi.S}.
\begin{proposition}\label{alg dep of int}
Let $E=F(\x_1,\x_2,\cdots,\x_t)$ be an antiderivative extension of $F$. An element $\x\in E$ is an antiderivative of $F$
if and only if there are a tuple $(c_1,\cdots,c_t)\in C^t$ and an element $f\in F$ such that $\x=\sum^t_{i=1}c_i\x_i+f.$
\end{proposition}

\bt \label{constructanti} Let $F$ be a differential field and suppose that there are elements $f_1,f_2,\cdots,f_n\in F$ satisfying the condition \textbf{C}. Let $E=F(\x_1,\x_2,\cdots,\x_n)$ be the field of rational functions over $F$ in variables $\x_1,\x_2,\cdots,\x_n$. Extend the derivation of $F$ to $E$ by defining $\x'_i=f_i$. Then $E$ is a no new constants extension of $F$. \et

\bpf Let $F_0:=F$ and $F_i:=F_{i-1}(\x_i)$. Suppose that the theorem is false. Then pick the smallest $k$ such that $F_k$ has a new constant. Note that $F_k=F_{k-1}(\x_k)$ and that there is an $a\in F_{k}-F_{k-1}$ such that $a'=0$. Therefore, by theorem \ref{new constant}, there is a $y\in F_{k-1}$ such that $y'=f_k\in F$. Now we apply proposition \ref{alg dep of int} and obtain that $y=\sum^{k-1}_{i=1}\a_i\x_i+f.$ Taking derivatives we obtain $\sum^k_{i=1}c_i f_i=f'$, where $c_1:=1$ and $c_i=-\a_i$ for all $i\ge 1$. This contradicts the choice of $f_i$'s.\epf

\begin{remark}\label{cx works}Let $F=\C(z)$ be the ordinary differential field of rational functions in one complex variable $z$ with the derivation $d/dz$. For any rational function $f\in F$, $\frac{df}{dz}$ has no simple pole and thus  for any distinct complex numbers $\a_1,\a_2,\cdots,\a_n$ and for any choice of constants $c_1, c_2,\cdots,c_n$, not all zero, there is no rational function $f\in F$ such that $\frac{df}{dz}=\sum^n_{i=1}c_i/(z+\a_i)$. Therefore, for the differential field $\C(z)$, elements $1/(z+\a_1), 1/(z+\a_2),\cdots,1/(z+\a_n)$ satisfy condition \textbf{C}.\end{remark}

\section{Picard-Vessiot Extensions for $\U(3,C)$}\label{lower unipotent groups}

In this section, we will provide a construction of a Picard-Vessiot extension for the group $\U(3,C)$.

\bprop\label{int a poly}
Let $K$ be a differential field with a field of constants $C$ and let $K(\x)$ be a no new constants extension of $K$ such that $\x$ is transcendental over $K$ and $\x'\in K$. Let $S=\sum^s_{i=0}q_i\x^i\in K[\x]$, $q_i\in K$ and $q_s\neq 0$.
\begin{itemize}

\item[1.] If there is a $T\in K(\x)$ such that $T'=S$ then $T\in K[\x]$ and its degree, deg $T$, equals $s$ or $s+1$.

\item[2.] And if $c\x'+f'\neq q_s$ for any $c\in C$ and for any $f\in K$ then there is no $T\in K(\x)$ such that $T'=S$.
\end{itemize}
\eprop

\bpf

Let there be an element $T\in K(\x)$ such that $T'=S$. Then there are relatively prime polynomials $P,Q\in K[\x]$, where $Q$ is monic, such that $T=P/Q$. Taking derivatives, we obtain \begin{equation}\label{no deno}Q^2S=P'Q-Q'P.\end{equation}
From the above equation, it is immediate that $Q$ divides $Q'$. On the other hand, since $Q$ is monic and $\x'\in K$ we know that deg $Q'<$ deg $Q$. Therefore  $Q=1$ and thus $P=T\in K[\x]$. Let $T=\sum^t_{i=0}r_i\x^i$, $r_t\neq 0$. From equation \ref{no deno},  we have \begin{equation}\label{integrability} r'_t\x^t+(tr_t\x'+r'_{t-1})\x^{t-1}+\cdots+r_1\x'+r'_0= \sum^s_{i=0}q_i\x^i.\end{equation}

Since $q_s\neq 0$, $t$ cannot be smaller than $s$. If $t\geq s+2$ then $r'_t=0$  and $tr_t\x'+r'_{t-1}=0$. Then $tr_t\x+r_{t-1}\in C\subset K$, contradicting the fact that $\x$ transcendental over $K$. Thus $t=s$ or $t=s+1$.

Furthermore, if $t=s$ then $r'_t=q_s$, where $r_t\in K$. And if $t=s+1$ then $r'_t=0$ and $tr_t\x'+r'_{t-1}=q_s$. Thus, we have shown that if
$ c\x'+f'\neq q_s$ for any $c\in C$ and  for any $f\in K$ then there is no $T\in K(\x)$ such that $T'=S$. \epf

\bt\label{generation step}
Let $F$ be a differential field with a field of constants $C$. Suppose that  $f_1,f_2\in F$ be elements satisfying the condition \textbf{C}. Let $E:=F(\x_1,\x_2)(\y)$ be the field of rational functions of $F$ in three variable $\x_1,\x_2$ and $\y$. Choose any $r\in F$ and extend the derivation of $F$  to  $E$ by setting \begin{equation}\label{3rd order der}\begin{pmatrix}
1 & \x_1 & \y \\
0 & 1 & \x_2\\
0 & 0 & 1
\end{pmatrix}'=\begin{pmatrix}
0 & f_1 & r \\
0 & 0 & f_2\\
0 & 0 & 0
\end{pmatrix}\begin{pmatrix}
1 & \x_1 & \y \\
0 & 1 & \x_2\\
0 & 0 & 1
\end{pmatrix}.\end{equation}
Then $E$ is a no new constants extension of $F$. In particular, there is no $y\in F(\x_1,\x_2)$ such that  $y'=\y'=f_1\x_2+r$ for any $r\in F$.
\et

\bpf
Suppose that the theorem is not true for some $r\in F$. Since $\x'_1=f_1$ and $\x'_2=f_2$, from theorem \ref{constructanti}, we know that $F(\x_1,\x_2)$ is a no new constants extension of $F$. Thus there is a new constant in the set $E-F(\x_1,\x_2)$. Note that $\y'=f_1\x_2+r\in F(\x_2)\subset F(\x_1,\x_2)$ and that $E=F(\x_1,\x_2)(\y)$. Therefore, by theorem \ref{new constant}, there is a $\z\in F(\x_1,\x_2)$ such that $\z'=f_1\x_2+r$.

Now, we apply proposition \ref{int a poly} (with $K=F(\x_1)$) and obtain that $\z\in F(\x_1)[\x_2]$ and that deg $\z$=1 or 2.

 Case deg $\z=2$:   Let $\z=r_2\x^2_2+r_1\x_2+r_0\in F(\x_1)[\x_2]$ with $r_2\neq 0$. Then taking derivatives, we obtain $r'_2=0$ and $2r_2f_2+r'_1=f_1$. This contradicts the condition \textbf{C}.

Case deg $\z=1$: Let $\z=r_1\x_2+r_0$. Then we have $r'_1=f_1$ and $r_1f_2+r'_0=r$. Thus there is a constant $c\in C$ such that $r_1=\x_1+c$ and therefore $(\x_1+c)f_2+r'_0=r$.

       Now we have \begin{equation}\label{imp eqn}f_2\x_1+s=-r'_0,\end{equation} where $r_0\in F(\x_1)$ and $s=cf_2-r\in F$. Then by proposition \ref{int a poly} (2), there are a constant $c\in C$ and an element $f\in F$ such that $cf_1+f'=f_2$. This again contradicts the condition \textbf{C}. \epf

\begin{corollary}\label{generation step cor}

Let $F$ and $E$ be as in theorem \ref{generation step} and let $\x\in E$ be an antiderivative of $F$. Then there are constants $c_1,c_2\in C$ and an element $s\in F$ such that $\x=c_1\x_1+c_2\x_2+s$.

\end{corollary}

\bpf
By proposition \ref{alg dep of int}, there is a constant $c\in C$ such that $\x+c\y\in F(\x_1,\x_2)$. We claim that $c=0$. Suppose not. Then, $c^{-1}\x+\y\in F(\x_1,\x_2)$ and let $s:=c^{-1}\x'+r\in F$. Let $E^*:=F(\x_1,\x_2)(\mu)$ be the field of rational functions in one variable $\z$ and extend the derivation of $F(\x_1,\x_2)$ to the field $E^*$ by setting $\z'=f_1\x_2+s$. But since $c^{-1}\x+\y\in F(\x_1,\x_2)$ and $(c^{-1}\x+\y)'=f_1\x_2+s$, we obtain that $E^*$ has a new constant, namely, $\z-(c^{-1}\x+\y)$. This contradicts theorem \ref{generation step}.

Thus $\x\in F(\x_1,\x_2)$. Now we again apply proposition \ref{alg dep of int} to prove the corollary.  \epf
In theorem \ref{heisen}, we will prove that the differential field $E$, as described in theorem \ref{generation step}, is a Picard-Vessiot extension with a differential Galois group  isomorphic to $\U(3,C)$ as groups. In order to do so, we will require the following two propositions to prove theorem \ref{heisen}.

\begin{proposition}\label{fixed field}
Let $F$ be a field and $E=F(\x_{i,j}|1\leq i\leq n, 1\leq j\leq n+1-i)$ be the field of rational functions over $n(n+1)/2$ variables.
Let $g$ be as in theorem \ref{Main result}. For $M:=\begin{pmatrix}
1 & c_{1,1} & c_{2,1}&\cdots& c_{n,1}\\
0 & 1 &  c_{1,2}& \cdots& c_{n-1,2}\\
\vdots & \vdots &  \vdots&\vdots& \vdots\\
0 & 0 &  \cdots&1 & c_{1,n}\\
0 & 0 &  \cdots & 0& 1
\end{pmatrix} \in \U(n+1,F)$ let $G:=\{\s_M:E\to E|M\in \U(n+1,F)\}$ be a collection of automorphisms on $E$ defined by $\s_M(g)=gM$. That is,
\begin{equation}\label{sigma defn}\s(\x_{i,j})=\left\{
  \begin{array}{ll}
   \x_{i,j}+\left(\sum^{i-1}_{t=1}c_{t,i+j-t}\ \x_{i-t,j}\right)+c_{i,j} , & \hbox{if \ $i\geq 2$;} \\
    \x_{1,j}+c_{1,j}, & \hbox{if \ $i=1$.}
  \end{array}
\right.\end{equation}

 Then $E^G$, the field fixed by $G$, equals $F$.
\end{proposition}

  \bpf
Let $S_{p,q}:=\{\x_{i,j}|1\leq i\leq p-1\}$
$\cup \{\x_{p,j}1\leq j\leq q-1\}$ and let $K_{p,q}:=F(S_{p,q})$.
Let $u\in E-F$ and choose the largest integer $p$ and a largest integer $q$ such that $u\in K_{p,q}(\x_{p,q})$. Then, there are relatively prime polynomials $P, Q\in K_{p,q}[\x_{p,q}]$ such that $u=P/Q$. Consider a matrix $M\in \U(n+1,F)$  such that $c_{i,j}=0$ for all $i\neq p$ and $j\neq q$, and $c_{p,q}\neq 0$.  Then  we see that $\s_{M}(\x_{p,q})=\x_{p,q}+c_{p,q}$ and  $\s_M(\x_{i,j})=\x_{i,j}$ for all $1\leq i\leq p$ and $j\leq q-1$.
Thus, in particular, the field $K_{p,q}$ is fixed by the automorphism $\s_M$.

Suppose that $\s_M(u)=u$. Then $\s_M(P)Q=\s_{M}(Q)P$ and since $P, Q$ are relatively prime, there is an element $r_{\s_{M}}\in K$ such that $\s_M(P(\x_{p,q}))=r_{\s_{M}} P(\x_{p,q})$ and $\s_M(Q(\x_{p,q}))=r_{\s_{M}} Q(\x_{p,q})$. That is $P(\x_{p,q}+c_{p,q})=r_{\s_{M}} P(\x_{p,q})$ and $Q(\x_{p,q}+c_{p,q})=r_{\s_{M}} Q(\x_{p,q})$. But these equations hold only when $c_{p,q}=0$, a contradiction.\epf

\begin{proposition}\label{coeff diffeq}
Let $E$ be a differential field extension of $F$.  Let $y_1,y_2,$ $\cdots,y_n$ $\in E$  and let $V:=$Span$_{C_E}\{y_1,y_2,\cdots,y_n\}$. Suppose that the wronskian $w(y_1,y_2,$ $\cdots,y_n)\neq 0$ and that the group $G:=G(E|F)$ of all differential automorphisms of $E$ fixing $F$ stabilizes the  vector space $V$. Then the differential operator $L(Y):=\dfrac{w(Y,y_1,y_2,\cdots,y_n)}{w(y_1,y_2,\cdots,y_n)}$ has coefficients in the differential field $E^{G}$. Moreover, $V$ is the full set of solutions of the differential equation $L(Y)=0$.
\end{proposition}

\bpf Write $L(Y)=Y^{(n)}+a_{n-1}Y^{(n-1)}+\cdots+a_1Y^{(1)}+a_0Y$. We will show that $a_i\in E^G$ for each $i$.  Clearly, $V$ is the full set of solutions of $L(Y)=0$. Let $L_{\s}:=Y^{(n)}+\s(a_{n-1})Y^{(n-1)}+\cdots+\s(a_1)Y^{(1)}+\s(a_0)Y$ for $\s\in G$. Since $G$ stabilizes $V$ and that $V$ is finite dimensional, we obtain that $G$ consists of automorphisms of the vector space $V$. It follows that $ker(L)=ker(L_{\s})=V$, in particular $ker(L-L_\s)\supset V$, and thus the dimension of $ker(L-L_\s)\geq n$. Since $(L-L_{\s})(Y)=(\s(a_{n-1})-a_{n-1})Y^{(n-1)}+\cdots+(\s(a_1)-a_1)Y^{(1)}+(\s(a_0)-a_0)Y$ is of order $\leq n-1$, we should have $L-L_{\s}=0$. That is, $a_i\in E^G$ for all $i, 0\leq i\leq n-1$.   \epf


\bt\label{heisen}

Let $E$ and $F$ be differential fields as defined in theorem \ref{generation step}.
Then $E$ is a Picard-Vessiot extension of $F$ for the differential operator \begin{equation}\label{3 order eqn}L(Y):=\frac{w(Y,1,\x_1,\y)}{w(1,\x_1,\y)},\end{equation} and $L^{-1}(0)=V$.  And, the differential Galois group $G:=G(E|F)$ is naturally isomorphic to the group $\U(3,C)$.
\et

\bpf
From theorem \ref{generation step}, we know that $E$ is a no new constants extension of $F$. Let $R:=F[\x_1,\x_2,\y]$ and let $\s\in G:=G(E|F)$. Since $\s(\x_i)'=\s(\x'_i)=\x'_i$, we have $\s(\x_1)=\x_1+\a_\s$ and  $\s(\x_2)=\x_2+\b_\s$, where $\a_\s,\b_\s\in C$. And since $\s(\y)'=\s(\y')=\s(\x'_1)\s(\x_2)+\s(f)=\x'_1(\x_2+\b_\s)+f=(\y+\b_\s\x_1)'$, there is a $\r_\s\in C$ such that $\s(\y)=\y+\b_\s\x_1+\r_\s$. Thus if $V:=$Span$_C\{1,\x_1,\y\}$ then $GV\subseteq V$. Let $g:=\begin{pmatrix}
1 & \x_1 & \y \\
0 & 1 & \x_2\\
0 & 0 & 1
\end{pmatrix}$ and $M_\s:=\begin{pmatrix}
1 & \a_\s & \r_\s \\
0 & 1 & \b_\s\\
0 & 0 & 1
\end{pmatrix}.$ Then we have $\s(g)=g M_\s$ and with respect to the basis $\{1,
\x_1,\y\}$, we have a group representation $\Gamma: G\to GL(3, C)$ defined by $\Gamma(\s)= M_\s$. Note that $F\langle V \rangle$, the differential field generated by $F$ and $V$, equals $E$. Therefore the action of $\s$ on the elements $\x_1$ and $\y$ completely determines the differential automorphism $\s$ on $E$. Thus $\Gamma$ is a faithfull representation of groups. Given a matrix $M:=\begin{pmatrix}
1 & \a & \r \\
0 & 1 & \b\\
0 & 0 & 1
\end{pmatrix}\in \U(3,C)$, we define a $F-$algebra automorphism $\s_M:R\to R$ such that $\s_M(g)=g M$.  Since $\s_M(\x_i)'=\s_M(\x'_i)$ for $i=1,2$ and $\s_M(\y)'=\s_M(\y')$, we see that $\s_M$ is an $F-$ algebra differential automorphism of $R$. Now one extends $\s_M$ to $E$, the field of fractions of $R$, to obtain a differential field automorphism of $E$. Hence $\Gamma (G)$ is isomorphic to $\U(3,C)$.

We note that the differential operator $L(Y):=\dfrac{w(Y,1,\x_1,\y)}{w(1,\x_1,\y)}$ has  span$_C$ $\{1,\x_1,$ $\y\}$ as the full set of solutions.  From proposition \ref{fixed field} and \ref{coeff diffeq}, we know that the coefficients of  the differential operator $L(Y)$  lie  in the field $F$.  Thus $E$ is a Picard-Vessiot extension of $F$ with Galois group isomorphic to $\U(3,C)$.    \epf

\begin{example}\label{heisen ex}Let $F=\C(z)$ with the derivation $d/dz$ and
 let $\a_1,\a_2$ be two distinct complex numbers. Let $f_i=1/(z+\a_i)$ for $i=1,2$, $r=0$ and let $E=F(\x_1,\x_2,\y)$ be the field of rational functions in three variables $\x_1,\x_2,\y$. Extend the derivation of $F$ to  $E$ using equation \ref{3rd order der}. Now we apply theorem \ref{heisen} and obtain that $E$ is a Picard-Vessiot extension of $F$ for the differential operator
\begin{equation}L:= \frac{d^3}{dz^3}+\frac{3z+\a_1+2\a_2}{(z+\a_1)(z+\a_2)}\frac{d^2}{dz^2}+\frac{1}{(z+\a_1)(z+\a_2)}\frac{d}{dz}, \end{equation} whose solution space $L^{-1}(0)=$Span$_\C\{1,\x_1,\y\}$, see \cite{Ravi.S w}. The differential Galois group of $E$ is isomorphic to the Heisenberg group $\U(3,C)$.  One can think of $\x_1,\x_2$ and $\y$ as $\log(z+\a_1)$, $\log(z+\a_2) $ and $\int^z_{a_0}\frac{\log(t+\a_2)}{t+\a_1}dt$, $a_0\neq \a_1$ respectively. The integral  $\int^z_{a_0}\frac{\log(t+\a_2)}{t+\a_1}dt$  is a `shifted' dilogarithm.
\end{example}


\section{Picard-Vessiot Extensions for $\U(n,C)$}

\bt\label{no new constants}
Let $F$ and $E$ differential fields as in theorem \ref{Main result}. Then
\begin{enumerate} \item [a.]\label{state a} The differential field $F\langle \x_{i,j}\rangle$ equals the field $F(\x_{1,i+j-1}, \x_{2, i+j-2},\cdots,$ $\x_{i-1,j+1}, \x_{i,j})$. In particular, $E=F\langle \x_{1,1},$ $\x_{2,1},\cdots, \x_{n,1}\rangle$.\\

\item[b.]\label{state b}  The differential field $N:=F\langle \x_{n,1}\rangle(\x_{1,1},\cdots,\x_{1,n})$ is a no new constants extension of $F$.\\

\item[c.]\label{state c} $E$ is a no new constants extension of $F$.\end{enumerate}\et
\bpf

 From equation \ref{defn derv}, we obtain \begin{align}
\label{UPVE deriv eqn1}\x'_{1j}&=f_j\\
\label{UPVE deriv eqn2}\x'_{ij}&=f_j\ \x_{i-1,j+1},\ \  2\leq i\leq n\ \text{and}\ 1\leq j\leq n+1-i.\end{align}

(a): Therefore $f_j\ \x_{i-1,j+1}=\x'_{i,j}\in F\langle \x_{i,j}\rangle$. And since $f_{j}\in F$, we obtain  $\x_{i-1,j+1}\in F\langle\x_{i,j}\rangle$. Repeating this argument, one proves that $F\langle \x_{i,j}\rangle\supseteq F(\x_{1,i+j-1},$ $\x_{2, i+j-2}, \cdots,\x_{i-1,j+1}, \x_{i,j})$. On the other hand, equations \ref{UPVE deriv eqn1} and \ref{UPVE deriv eqn2} also tell us that $F(\x_{1,i+j-1}, \x_{2, i+j-2},$ $ \cdots,\x_{i-1,j+1}, \x_{i,j})$ is a differential field. Since $F\langle\x_{i,j}\rangle$ is the smallest differential field containing $F$ and $\x_{i,j}$, we obtain $F\langle \x_{i,j}\rangle$ $=F(\x_{1,i+j-1}, \x_{2, i+j-2},$ $\cdots,\x_{i-1,j+1},$  $\x_{i,j})$. It is easy to check that $E=F\langle \x_{1,1},$ $\x_{2,1},\cdots, \x_{n,1}\rangle$.

(b): Let $N_k:= F(\x_{1,n-k},\cdots,\x_{1,n})(\x_{2,n-1},\cdots,\x_{k+1,n-k})$ and observe from statement (a) that $N_k=F(\x_{1,n-k},\cdots,\x_{1,n})\langle\x_{k+1,n-k)}\rangle$.
We see from theorem \ref{generation step} that $N_1$ is a no new constants extension of $F$. Assume that $N_k$  is a no new constants extension of $F$ for some $k\geq 1$. Let $K=F(\x_{1,n-k},\cdots,\x_{1,n})$ $(\x_{2,n-1},\cdots,\x_{k,n-(k-1)})$ and note that $N_{k+1}=K(\x_{1,n-(k+1)}$, $\x_{k+1,n-k})$ $(\x_{k+2,n-(k+1)})$. Applying theorem \ref{generation step} we obtain $N_{k+1}$ is a no new constants extension of $K$. Since $K\subset N_{K}$, we see that $K$ is a no new constants extension of $F$. Thus $N_k$ is a no new constants extension of $F$. Choose $k=n$ to prove statement (b).

(c): The case $n=2$ follows from theorem \ref{generation step}. Assume that (c) is true for some $n\geq 3$  and let $F^*:=F\langle \x_{n,1}\rangle$, \begin{equation*} g^*:=\begin{pmatrix}
1 & \x_{1,1} & \x_{2,1}&\cdots& \x_{n-1,1}\\
0 & 1 &  \x_{1,2}& \cdots& \x_{n-1,2}\\
\vdots & \vdots &  \vdots&\vdots& \vdots\\
0 & 0 &  \cdots&1 & \x_{1,n}\\
0 & 0 &  \cdots & 0& 1
\end{pmatrix},\ \text{and} \ \ A^*:=\begin{pmatrix}
0 & f_1 & 0&\cdots& 0\\
0 & 0 &  f_2& \cdots& 0\\
\vdots & \vdots &  \vdots&\vdots& \vdots\\
0 & 0 &  \cdots&0 & f_{n-1}\\
0 & 0 &  \cdots & 0& 0
\end{pmatrix}.\end{equation*}
Then from (a) it follows that  $E=F^*(g^*)$. Since ($g^*)'=A^*g^*$, applying induction, we obtain that $E$ is a no
new constants extension of $F^*$. From (b) we know that $N$ is a no new constants extension of $F$. In particular,
$F^*\subset N$ is a no new constants extension of $F$ as well.  Thus we have shown that $E$ is a no new constants extension of $F$.\epf
\begin{corollary}\label{nnc induction cor}
Let $\x\in N$ be an antiderivative of $F$. Then, there are constants $c_1,\cdots,c_n\in C_F$ and an element $f\in F$ such that $\x=\sum^n_{i=1}c_i\x_i+f$. Moreover, if $\x\in F\langle \x_{1,n}\rangle$ then $c_{i}=0$ for all $i\leq n-1$.

\end{corollary}

\bpf

From corollary \ref{generation step cor} it is enough to consider the case when $n\geq 3$. Assume the corollary for the field $N^*:=F(\x_{1,2},\cdots,\x_{1,n})\langle \x_{n-1,2}\rangle$. Let $K:=F(\x_{1,2},\cdots,\x_{1,n})\langle \x_{n-2,3}\rangle$ and note that $K\subset N^*$. Since $\x\in N=K(\x_{1,1},\x_{n-1,2})(\x_{n,1})$ and $\x'\in F\subset K$, by corollary \ref{generation step cor}, there are constants $c_1,d_1\in C$ such that $\x-c_1\x_{1,1}-d_1\x_{n-1,2}\in K\subset N^*$.  Since $\x_{n-1,2}\in N^*$, we have $\x-c_1\x_{1,1}\in N^*$. Note that $\x-c_1\x_{1,1}\in N^*$ is an antiderivative of $F$ and therefore, applying induction to the differential field $N^*$, we obtain constants $c_2,\cdots,c_n$ such that $\x-c_1\x_{1,1}=\sum^n_{i=2}c_i\x_{1,i}+f$ for some $f\in F$ as desired. From (a), we know that $F\langle \x_{1,n}\rangle=F(\x_{1,n},\x_{2,n-1},\cdots,\x_{n,1})$ and thus  $\x_{1,1},\cdots,\x_{1,n-1}$ remains algebraically independent over $F\langle \x_{1,n}\rangle$. Hence $\x=\sum^n_{i=1}c_i\x_{1,i}+f\in F\langle \x_{1,n}\rangle$ implies $c_{i}=0$ for all $i\leq n-1$.  \epf


\begin{corollary}\label{no new constants cor}
Let $\x\in E$ be an antiderivative of $F$. Then there are constants $c_1,\cdots,c_n\in C$ and an $f\in F$ such that $\x=\sum^n_{i=1}c_i\x_{i}+f$. Moreover, if $\x\in F\langle \x_{s,t}\rangle$ then $c_i=0$ for all $i\neq s+t-1$ and thus $\x=c_{s+t-1}\x_{1,s+t-1}+f$.

\end{corollary}

\bpf

Since $\x\in E=F^* (g^*)$, applying induction with $F\langle \x_{n,1}\rangle$ as our base field, we obtain constants $c_1,\cdots,c_{n-1}\in C$ and an $\tilde{f}\in F\langle \x_{n,1}\rangle$ such that \begin{equation}\label{no new constants cor eqn}\x=\sum^{n-1}_{i=1}c_i\x_{1,i}+\tilde{f}.\end{equation} Then $\tilde{f}=\x-\sum^{n-1}_{i=1}c_i\x_{1,i}\in F\langle \x_{n,1}\rangle$ is an antiderivative of $F$. Now we apply corollary \ref{nnc induction cor} to the differential field $F\langle \x_{n,1}\rangle$ and obtain a constant $c_n\in C$ and an element $f\in F$ such that $\tilde{f}=c_n\x_{1,n}+f$. Substituting back for $\tilde{f}$ in equation \ref{no new constants cor eqn}, we obtain $\x=\sum^n_{i=1}c_i\x_{1,i}+f$.

Note that $\x_{1,1},\cdots,\x_{1,s+t-2},$ $\x_{1,s+t},\cdots,\x_{1,n}$ remains algebraically independent over $F\langle \x_{s,t}\rangle$, see (a). Thus $\x=\sum^n_{i=1}c_i\x_{1,i}+f\in F\langle \x_{s,t}\rangle$ implies that $c_i=0$ for all $i\neq s+t-1$.\epf



\textsl{Proof of theorem \ref{Main result}}.  We know that $E$ is a no new constants extension of $F$.
For $\s\in G$ and for each $1\leq j\leq n$, we have $\s(\x_{1,j})=\x_{1,j}+c^{\s}_{1,j}$ for some constants $c^{\s}_{1,j}\in C$. For $1\leq j\leq n-1$, we have $\x'_{2,j}= f_j\x_{1,j+1}$ and therefore $$\s(\x_{2,j})'=\s(\x'_{2,j})= f_j\ \s(\x_{1,j+1})= f_j\ \x_{1,j+1}+ f_j\ c^{\s}_{1,j+1}=(\x_{2,j}+c^{\s}_{1,j+1}\ \x_{1,j})'.$$
Thus, there are constants $c^{\s}_{2,j}\in C$ such that $\s(\x_{2,j})=\x_{2,j}+c^{\s}_{1,j+1}\ \x_{1,j}+c^{\s}_{2,j}$, for all $1\leq j\leq n-2$. Assume that for some integer  $s, 2\leq s\leq n$, there are constants $c^{\s}_{i,j}\in C$ such that $$\s(\x_{s,j})=\x_{s,j}+\left(\sum^{s-1}_{t=1}c^{\s}_{t,s+j-t}\ \x_{s-t,j}\right)+c^{\s}_{s,j},$$
for all $1\leq j\leq n+1-s.$
Note that \begin{align*}\s(\x_{s+1,j})'&= f_j\s(\x_{s,j+1})\\
&= f_j\left(\x_{s,j+1}+\left(\sum^{s-1}_{t=1}c^{\s}_{t,s+1+j-t}\ \x_{s-t,j+1}\right)+c^{\s}_{s,j+1}\right).\\
&=\left(\x_{s+1,j}+\left(\sum^{s-1}_{t=1}c^{\s}_{t,s+1+j-t}\ \x_{s+1-t,j}\right)+c^{\s}_{s,j+1}\ \x_{1,j}  \right)'
\end{align*}
and thus there is a $c^{\s}_{s+1,j}\in C$ such that \begin{align*}\s(\x_{s+1,j})&=\x_{s+1,j}+\left(\sum^{s-1}_{t=1}c^{\s}_{t,s+1+j-t}\ \x_{s+1-t,j}\right)+c^{\s}_{s,j+1}\ \x_{1,j}+c^{\s}_{s+1,j}\notag,\\
&=\x_{s+1,j}+\left(\sum^{s}_{t=1}c^{\s}_{t,s+1+j-t}\ \x_{s+1-t,j}\right)+c^{\s}_{s+1,j}.
\end{align*}

Then, by induction, for any fixed $i, 1\leq i\leq n$ and for all $j,1\leq j\leq n+1-i$ there are constants $c^{\s}_{i,j}\in C$ such that
\begin{equation}\label{rep eqn}\s(\x_{i,j})=\left\{
  \begin{array}{ll}
   \x_{i,j}+\left(\sum^{i-1}_{t=1}c^{\s}_{t,i+j-t}\ \x_{i-t,j}\right)+c^{\s}_{i,j} , & \hbox{if \ $i\geq 2$;} \\
    \x_{1,j}+c^{\s}_{1,j}, & \hbox{if \ $i=1$.}
  \end{array}
\right.\end{equation}

Thus, if  $g:=\begin{pmatrix}
1 & \x_{1,1} & \x_{2,1}&\cdots& \x_{n,1}\\
0 & 1 &  \x_{1,2}& \cdots& \x_{n-1,2}\\
\vdots & \vdots &  \vdots&\vdots& \vdots\\
0 & 0 &  \cdots&1 & \x_{1,n}\\
0 & 0 &  \cdots & 0& 1
\end{pmatrix}$ and $\s\in G$ then, from equation \ref{rep eqn}, we see that there is an element $$M_\s:=\begin{pmatrix}
1 & c^\s_{1,1} & c^\s_{2,1}&\cdots& c^\s_{n,1}\\
0 & 1 &  c^\s_{1,2}& \cdots& c^\s_{n-1,2}\\
\vdots & \vdots &  \vdots&\vdots& \vdots\\
0 & 0 &  \cdots&1 & c^\s_{1,n}\\
0 & 0 &  0&\cdots& 1
\end{pmatrix} \in \U(n+1, C)$$ such that $\s(g)=g M_\s $.

Let $V:=$Span$_C\{1,\x_{1,1},\x_{2,1},\cdots,\x_{n,1}\}$. From equation \ref{rep eqn}, it is clear that $GV\subset V$. And, with respect to the basis $1,\x_{1,1},\x_{2,1},\cdots,\x_{n,1}$,  we have a  representation $\Gamma : G\to GL(n+1, C)$ defined by $\Gamma(\s)=M_\s$. And indeed,
 $\Gamma(G)\subseteq \U(n+1,C)$. Note that if $\s,\rho\in G$ agrees on $V$, then  they agree on $F\langle V\rangle=E$ and thus $\Gamma$ is injective. For any matrix \begin{equation}\label{surjection} M:=\begin{pmatrix}
1 & c_{1,1} & c_{2,1}&\cdots& c_{n,1}\\
0 & 1 &  c_{1,2}& \cdots& c_{n-1,2}\\
\vdots & \vdots &  \vdots&\vdots& \vdots\\
0 & 0 &  \cdots&1 & c_{1,n}\\
0 & 0 &  0&\cdots& 1
\end{pmatrix} \in \U(n+1,C)\end{equation} let $\s_M$ be the  $F-$algebra automorphism on $F[g]$ defined by $\s_M(g)=g M$. Note that $\s_M\s_{M^{-1}}=\s_{M^{-1}}\s_{M}=I$, the identity matrix. Thus $\s_M$ is a ring automorphism of $F[g]$. To show that $\s_M$ is a differential automorphism, we only need to check that $\s_M(\x_{i,j})'=\s_M(\x'_{i,j})$. Since $\s_M(\x_{1,j})=\x_{1,j}+c_{1,j}$, we obtain $\s_M(\x_{1,j})'= \x'_{1,j}= \s_M(\x'_{1,j})$. Now \begin{align*}\s_M(\x'_{i+1,j})&=\s_M( f_j\x_{i,j+1})\\
&= f_j\s_M(\x_{i,j+1})\\
 &= f_j\left(\x_{i,j+1}+\left(\sum^{i-1}_{t=1}c_{t,i+1+j-t}\ \x_{i-t,j+1}\right)+c_{i,j+1}\right)\\
&=\left(\x_{i+1,j}+\left(\sum^{i-1}_{t=1}c_{t,i+1+j-t}\ \x_{i+1-t,j}\right)+c_{i,j+1}\ \x_{1,j}  \right)'\\
&=\s_M(\x_{i+1,j})'.
 \end{align*}
 Thus $\s_M$ is a differential automorphism of the ring $F[g]$.
 Now extend $\s_M$ to a differential field automorphism of the field of fractions $E$ of $F[g]$. Thus $G$ is isomorphic to the group $\U(n+1, C)$. From propositions \ref{fixed field} and \ref{coeff diffeq}, we obtain that the differential operator $$L(Y):=w(Y,1,\x_{1,1},\x_{2,1},\cdots,\x_{n,1})/w(1,\x_{1,1},\x_{2,1},\cdots,\x_{n,1})$$ has coefficients in the field $F$. Since $V$ is the full set of solutions of the differential equation $L(Y)=0$ and $F\langle V\rangle=E$, we conclude that $E$ is a Picard-Vessiot extension of $F$ for the differential operator $L(Y)$. It is also clear that one can realize any full unipotent subgroup of $GL(n,C)$ by suitably choosing a basis for $V$.  $\square$

\bprop \label{pv ring}

The differential ring $F[g]$ is the Picard-Vessiot ring of $E$. 

\eprop

\bpf
It suffices to show that $F[g]$ is a simple differential ring. Let $I$ be a  differential ideal of $F[g]$ and suppose that $I\cap F=\{0\}$. Let $S_{p,q}:=\{\x_{i,j}|1\leq i\leq p-1\}$
$\cup \{\x_{p,j}|1\leq j\leq q-1\}$. It is easy check that $F[S_{p,q}]$ is a differential ring. Choose $p, q$ such that $I\cap F[S_{p,q}][\x_{p,q}]\neq \{0\}$ and that $I\cap F[S_{p,q}]=\{0\}$. Choose a non zero element $u\in I\cap F[S_{p,q}][\x_{p,q}]$ of smallest possible degree  and write $u=\sum^n_{i=0}a_i\x^i_{p,q}$, where $a_n\neq 0$, $a_i\in F[S_{p,q}]$. Note that $n\geq 1$ , degree of $a'_nu-a_nu' \leq n-1$ and that $a'_nu-a_nu'\in I$. Thus $a'_nu-a_nu'=0$ and in particular, $a'_na_{n-1}-a_n(na_nf_q\x_{p-1,q+1}+a'_{n-1})=0$. Then $(-a_{n-1}/na_n)'=f_q\x_{p-1,q+1}$ in $E$, which implies $\x_{p,q}+(a_{n-1}/na_n)$ is a new constant. This contradicts the fact that $E$ is a no new constants extension of $F$.\epf

\section{Examples}

\subsection{Hyperlogarithms}
Consider the differential field $F=\C(z)$ with the usual derivation $d/dz$. Let $\a_1,\cdots,\a_n$ be distinct complex numbers. A hyperlogarithmic function is an iterated integral of the form
\begin{equation}
\label{hyp form} L(\a_1,\a_2,\cdots,\a_n|z,z_0):=\int^z_{z_0}\int^{s_{n-1}}_{z_0}\cdots\int^{s_1}_{z_0}\frac{ds_0}{s_0-\a_n}\cdots \frac{ds_{n-1}}{s_{n-1}-\a_1}
\end{equation}
where $z_0$ is a fixed point and $z_0\neq \a_n$, see \cite{hyperlog}. Let $\mathcal{M}(U)$ be the field of meromorphic functions on $U$, where $z_0\in U$ is a simply connected domain that does not contain $\a_i$ for any $i$.
Let $f_j:=1/(z+\a_j)$. As noted in remark \ref{cx works},  these $f_j$'s satisfy the condition \textbf{C}.  Using theorem \ref{Main result} we may construct a Picard-vessiot extension $F(\x_{i,j}|1\leq i\leq n, 1\leq j\leq n+1-i)$ for the differential operator \begin{equation}\label{cx eqn}L(Y)=w(Y,1,\x_{1,1},\x_{2,1},\cdots,\x_{n,1})/w(1,\x_{1,1},\x_{2,1},\cdots,\x_{n,1}).\end{equation}

Define an $F$-algebra homomorphism $\phi:=F[g]\to \mathcal{M}(U)$ such that $\phi(\x_{p,q})=L(\a_q,\cdots, \a_{p+q-1}|z,z_0)$. Then, one can see that this map commutes with the derivation $d/dz$ and therefore ker $\phi$ is a differential ideal of $F[g]$. Applying proposition \ref{pv ring}, we conclude that ker $\phi=\{0\}$ and thus $\phi$ is injective. This shows that the collection
$$\{L(\a_q,\a_{q+1},\cdots,\a_{p+q-1}|z,z_0)| \  \ 1\leq i\leq p, 1\leq j\leq n+1-p \}$$
of hyperlogarithms is algebraically independent over $\C(z)$.



We may compute the differential equation \ref{cx eqn} using a Maple program, see \cite{Ravi.S w}. Here, I will list differential equations for groups $\U(n+1,\C),$ when $n=2,3$ and $4$.

 $$n=2\quad\quad\frac{d^3}{dz^3}+\frac{3z+\a_1+2\a_2}{(z+\a_1)(z+\a_2)}\frac{d^2}{dz^2}+\frac{d}{dz}.$$

\begin{align*}n=3\quad\quad &\frac{d^4}{dz^4}+\frac{6 z^2+(3\a_1+4 \a_2+5 \a_3)z+2 \a_3\a_1+\a_2 \a_1+3 \a_3\a_2}{(z+\a_1)(z+\a_2)(z+\a_3)}\frac{d^3}{dz^3}\\&+\frac{7z+\a_1+2\a_2+4\a_3}
{(z+\a_1)(z+\a_2)(z+\a_3)}\frac{d^2}{dz^2}+\frac{1}{(z+\a_1)(z+\a_2)(z+\a_3)}\frac{d}{dz}.\end{align*}

 \begin{align*}n=4\quad\quad\frac{d^5}{dz^5}&+\frac{P_1 z^3+P_2 z^2+ P_3 z+P_4}{S}\frac{d^4}{dz^4}+\frac{7z+\a_1+2\a_2+4\a_3}
{S}\frac{d^3}{dz^3}\\&+\frac{P_5 z^2+P_6z+P_7}{S}  \frac{d^2}{dz^2}+\frac{1}{S}\frac{d}{dz},\end{align*}
where $P_1=10$, $P_2=9\a_4+8 \a_3+7 \a_2+6\a_1$,$P_3=5\a_1\a_4+4\a_1\a_3+5\a_2\a_5+6\a_2\a_4+3\a_2\a_1+7\a_3a_4$, $P_4=\a_1\a_2 \a_3+2\a_1\a_2 \a_4+3 \a_1\a_3\a_4+4\a_2\a_3\a_4$, $P_5=25$, $P_6=7\a_1+10\a_2+14\a_3+19\a_4$, $P_7=\a_1\a_2+2\a_1\a_3+3\a_2\a_3+9\a_3\a_4+6\a_2\a_4+4a_1\a_4$ and $S=\prod^4_{i=1}(z+\a_i)$.




\begin{thebibliography}{99}


\bibitem{j.ax}
    J. Ax, {\em On Schanuel's Conjectures},
    Ann. of Math (2) \textbf{93} (1971), 252-268. MR \textbf{43}

\bibitem{Bi-Bi}

A. Bialynicki-Birula, {\em On the inverse problem of Galois theory of differential fields}, Bull.
Amer. Soc. \textbf{16} (1963), 960-964.


\bibitem{hyperlog} Matthieu Deneufch\^{a}tel , G\'{e}rard Henry Edmond Duchamp, Vincel Hoang Ngoc Minh, Allan I. Solomon, 
{\em Independence of hyperlogarithms over function fields via algebraic combinatorics}, arXiv:1101.4497v1 [math.CO].



\bibitem{Magid 1}
    A. Magid, \ {\bf Lectures on Differential Galois Theory},
  University Lecture Series. American Mathematical society
  1994, 2nd edn.

    \bibitem{Marius-Singer}
M. van der Put, M. F. Singer, {\bf Galois Theory of Linear Differential Equations}, \textbf{328}, Grundlehren der mathematischen Wissenshaften, Springer, Heidelberg, 2003.

\bibitem{mit-sing}

C. Mitschi, M. Singer, {\em Connected Linear Groups as Differential Galois Groups}, \text{184} (1996), 333-361.

\bibitem{M.Ros}
    M. Rosenlicht, {\em On Liouville's Theory of Elementary Functions},
     Pacific J. Math (2)  \textbf{65} (1976), 485-492.


 \bibitem{Rosenlicht-Singer}

M. Rosenlicht, M. Singer, {\em On Elementary, Generalized Elementary, and Liouvillian Extension Fields}, Contributions to Algebra, (H. Bass et.al., ed.), Academic Press (1977) 329-342 .


\bibitem{Ravi.S}

V. Ravi Srinivasan, {\em Iterated Antiderivative Extensions}, Journal of Algebra (8) \textbf{324} (2020) 2042-2051.

\bibitem{Ravi.S w}
V. Ravi Srinivasan, {\em \url{http://andromeda.rutgers.edu/~ravisri/}}.








\end{thebibliography}
\end{document}